\documentclass[12pt,francais,nothms]{aclart}
\usepackage{mathptm}
\usepackage[all]{xy}

\usepackage[all]{xy}

\RequirePackage[T1]{fontenc}
\RequirePackage{mathrsfs,amsmath,amssymb}
\let\mathcal\mathscr
 
\RequirePackage{url}

\RequirePackage{xspace}
\RequirePackage[english,frenchb]{babel}

\usepackage{amssymb}

\newtheorem{e-proposition}[theorem]{Proposition}

\newtheorem{e-definition}[theorem]{Definition\rm}

\newtheorem{def-lem}[theorem]{Lemma-Definition}
\newtheorem{def-prop}[theorem]{Proposition-Definition}

\newtheorem{theoreme}{Th\'eor\`eme}[section]

\newtheorem{proposition}[theoreme]{Proposition}

\newtheorem{fdef-lem}[theoreme]{Lemme-D\'efinition}
\newtheorem{fdef-prop}[theoreme]{Proposition-D\'efinition}
\renewcommand{\theequation}{\arabic{equation}}
\def\og{\leavevmode\raise.3ex\hbox{$\scriptscriptstyle\langle\!\langle$~}}
\def\fg{\leavevmode\raise.3ex\hbox{~$\!\scriptscriptstyle\,\rangle\!\rangle$}}
\def\ord{{\operatorname {\rm ord}}}
\def\mod{\ {\rm mod}\ }
\def\ac{{\overline{\rm ac}}}
\def\11{{\mathbf 1}}

\def\Def{{\rm Def}}

\def\RDef{{\rm RDef}}
\def\RDefe{{\rm RDef}^{\rm exp}}
\def\RDefw{{\rm RDef}^{\rm e}}

\def\LPas{\cL_{\rm DP}}

\def\LO{{\cL_{\cO}}}

\def\Fv{{\mathfrak{F}}}

\def\AA{{\mathbf A}}

\def\CC{{\mathbf C}}

\def\ee{{\mathbf e}}

\def\GG{{\mathbf G}}

\def\LL{{\mathbf L}}

\def\NN{{\mathbf N}}

\def\QQ{{\mathbf Q}}

\def\ZZ{{\mathbf Z}}

\def\cA{{\mathcal A}}
\def\cB{{\mathcal B}}
\def\cC{{\mathcal C}}
\def\cD{{\mathcal D}}

\def\cI{{\mathcal I}}

\def\cL{{\mathcal L}}

\def\cO{{\mathcal O}}
\def\cP{{\mathcal P}}

\def\cS{{\mathcal S}}

\def\llp{\mathopen{(\!(}}
\def\llb{\mathopen{[\![}}
\def\rrp{\mathopen{)\!)}}
\def\rrb{\mathopen{]\!]}}

\begin{document}





%
\selectlanguage{french}
\title[Fonctions constructibles  exponentielles]{Fonctions constructibles  exponentielles, transformation de Fourier motivique et principe de transfert}
\alttitle{Constructible exponential functions, motivic Fourier transformation and transfer principle}



\author{Raf Cluckers}
\address{{\'E}cole Normale Sup{\'e}rieure,
D{\'e}partement de math{\'e}matiques et applications, UMR 8553 du
CNRS, 45, rue d'Ulm, 75230~Paris cedex 05, France}
\email{raf.cluckers@ens.fr}
\author{Fran\c cois Loeser}
\address{{\'E}cole Normale Sup{\'e}rieure,
D{\'e}partement de math{\'e}matiques et applications, UMR 8553 du
CNRS, 45, rue d'Ulm, 75230~Paris cedex 05, France}
\email{Francois.Loeser@ens.fr}

\begin{altabstract}\selectlanguage{english}
We introduce spaces of exponential constructible functions in the motivic setting for which
we construct direct image functors  in the absolute and relative cases.
This allows us to define a motivic
Fourier transformation for which we get various inversion statements.
We define also motivic
Schwartz-Bruhat spaces on which motivic
Fourier transformation induces an  isomorphism.
Our motivic integrals specialize to non archimedian integrals.
We give a general transfer principle 
comparing identities between functions defined by integrals over local fields
of characteristic zero,  resp. positive,  having the same residue field.
Details of constructions and proofs will
be given elsewhere.
\end{altabstract}

\begin{abstract}
\selectlanguage{french}
Nous d\'efinissons des espaces de fonctions constructibles exponentielles dans le cadre motivique
pour lesquels nous construisons des foncteurs d'image directe dans le cas absolu et relatif.
Ceci nous permet de d\'efinir une  transformation de Fourier motivique pour laquelle nous obtenons
des th\'eor\`emes d'inversion. Nous d\'efinissons \'egalement des espaces de Schwartz-Bruhat motiviques
sur lesquels la transformation de Fourier motivique induit un isomorphisme.
Nos int\'egrales motiviques se sp\'ecialisent sur des int\'egrales non archim\'ediennes.
On donne un principe g\'en\'eral de transfert 
comparant les  identit\'es entre fonctions d\'efinies par des int\'egrales
sur des corps locaux de caract\'eristique z\'ero ou positive ayant m\^eme corps r\'esiduel.
Les d\'etails des constructions et des preuves seront donn\'es
ailleurs.
\end{abstract}

\maketitle

\selectlanguage{english}
\section*{Abridged English version}We keep notations and framework from  \cite{cr1},  \cite{cr2}
and \cite{cl}.
In particular, we 
fix a field $k$ of characteristic zero and we consider fields $K$ containing $k$ and the field of Laurent series $K \llp t \rrp$
endowed with its standard valuation  $\ord : K \llp t \rrp^{\times} \rightarrow \ZZ$.
For $x$ in  $K \llp t \rrp$ we set  $\ac (x) = x t^{-\ord (x)} \mod t$ if $x \not=0$
and $\ac (0) = 0$, which gives $\ac : K \llp t \rrp \rightarrow K$.
We denote by  $\Def_k$ the category of  definable $k$-subassignments
in the  Denef-Pas language $\LPas$. 
For every object  $Z$ we defined in  loc. cit. categories $\Def_{Z}$ and
$\RDef_{Z}$ consisting of definable subassignments over  $Z$.

Let $Z$ be in $\Def_k$. We consider the category $\RDefe_{Z}$
whose objects are triples $(Y\to Z, \xi, g)$ with $Y$ in
$\RDef_{Z}$ and $\xi : Y \rightarrow h[0,1,0]$ and $g : Y
\rightarrow h[1,0,0]$ morphisms in $\Def_k$. A morphism
$(Y'\to Z, \xi', g') \rightarrow (Y\to Z, \xi, g)$
in $\RDefe_{Z}$ is a morphism $h : Y' \rightarrow Y$ in $\Def_Z$ such that
$\xi' = \xi \circ h$ and $g' = g \circ h$.
The functor sending 
$Y$ in $\RDef_{Z}$ to $(Y, 0, 0)$, with $0$ denoting the constant morphism
with value $0$ in $h[0,1,0]$, resp. $h[1, 0,0]$
being fully faithful, we may consider $\RDef_{Z}$ as a full subcategory of
$\RDefe_{Z}$. We shall also consider
the intermediate 
full subcategory $\RDefw_{Z}$ consisting of objects 
$(Y, \xi, 0)$ with $\xi : Y \rightarrow h[0,1,0]$ a morphism in $\Def_k$.
To the category $\RDefe_{Z}$ one assigns a Grothendieck ring
$K_0 (\RDefe_{Z})$ defined as follows. As an abelian group it is 
 the
quotient of the free abelian group over symbols $[Y \rightarrow Z,
\xi,g]$ with $(Y \rightarrow Z, \xi,g)$ in $\RDefe_{Z}$ by the
following 
four 
relations
\begin{equation}[Y \rightarrow Z, \xi,g] = [Y' \rightarrow Z, \xi',g']
\end{equation}
for $(Y \rightarrow Z, \xi,g)$ isomorphic to $(Y' \rightarrow Z,
\xi',g')$,
\begin{equation}
[(Y \cup Y') \rightarrow Z, \xi, g]
+ [(Y \cap Y') \rightarrow Z,
\xi_{|Y \cap Y'},g_{|Y \cap Y'}] 
= [Y \rightarrow Z, \xi_{|Y},
g_{|Y}] + [Y' \rightarrow Z, \xi_{|Y'}, g_{|Y'}]
\end{equation}
for $Y$ and $Y'$ definable subassignments of some $X$ in $\RDef_Z$
and $\xi$, $g$ defined on $Y \cup Y'$,
 \begin{equation} [Y\to Z,\xi,g+h]=[Y\to Z,\xi + \overline h,g]
 \end{equation}
for $h:Y\to h[1,0,0]$ a definable morphism with $\ord(h(y)) \geq 0$ for
all $y$ in $Y$ and $\overline h$ the reduction of $h$ modulo $t$,
and
 \begin{equation}
 [Y[0,1,0]\to Z,\xi+p,g]=0
 \end{equation}
when $p:Y[0,1,0]\to h[0,1,0]$ is the projection and when
$Y[0,1,0]\to Z$, $g$, and $\xi$ factorize through the projection
$Y[0,1,0]\to Y$.
Fiber product endows $K_0 (\RDefe_{Z})$ with a ring structure.

Finaly, one defines the ring 
$\cC  (Z)^{\rm exp}$ of exponential constructible functions
as
$\cC  (Z)^{\rm exp}:=\cC (Z)\otimes_{K_0 (\RDef_{Z})} K_0 (\RDefe_{Z})$.
One defines similarly $C (Z)^{\rm exp}$ and ${\rm I}_S C (Z)^{\rm exp}$.
In our main result, Theorem \ref{mte},  we construct direct image morphisms
$f_! : {\rm I}_S C (Z)^{\rm exp} \rightarrow
{\rm I}_S C (Y)^{\rm exp}$ for $f:  Z \rightarrow Y$ in $\Def_S$, extending those constructed in loc. cit. 
This is extended to the relative setting in Theorem \ref{mter}.
This formalism allows us to define in section \ref{fourier} a motivic Fourier transformation 
for which we can prove several inversion Theorems. In particular we define motivic Schwartz-Bruhat
spaces on which Fourier transformation induces an automorphism.
Our motivic integrals specialize to non archimedian integrals.
We give a general transfer principle comparing identities between functions defined by integrals over local fields
of characteristic zero,  resp. positive characteritic, having the same residue field.


\setcounter{equation}{0}
\selectlanguage{french}

\section{Fonctions constructibles exponentielles}

\subsection{Anneaux de Grothendieck exponentiels}\label{grthe}On reprend
le cadre et les notations des notes \cite{cr1} 
et \cite{cr2}
et de l'article \cite{cl}.
En particulier on
fixe un corps $k$ de caract\'eristique
z\'ero et on consid\`ere des corps $K$ contenant $k$ ainsi que le
corps des s\'eries de Laurent $K \llp t \rrp$
muni
de la valuation naturelle  $\ord : K \llp t \rrp^{\times} \rightarrow \ZZ$.
Pour $x$ dans $K \llp t \rrp$ on pose $\ac (x) = x t^{-\ord (x)} \mod t$ si $x \not=0$
et $\ac (0) = 0$, ce qui donne $\ac : K \llp t \rrp \rightarrow K$.
On note $\Def_k$ la cat\'egorie des $k$-sous-assignements d\'efinissables
dans  le langage de Denef-Pas $\LPas$. 
Pour chaque objet $Z$ on a d\'efini dans loc. cit. des cat\'egories $\Def_{Z}$ et
$\RDef_{Z}$ form\'ees de sous-assignements d\'efinissables au-dessus de $Z$.

On consid\`ere maintenant la cat\'egorie $\RDefe_{Z}$
dont les objets sont les triplets $(Y\to Z, \xi, g)$ avec  $Y$ dans
$\RDef_{Z}$, $\xi : Y \rightarrow h[0,1,0]$ et $g : Y
\rightarrow h[1,0,0]$ des morphismes dans  $\Def_k$. Un morphisme
$(Y'\to Z, \xi', g') \rightarrow (Y\to Z, \xi, g)$
dans $\RDefe_{Z}$ est un  morphisme $h : Y' \rightarrow Y$ dans $\Def_Z$ tel que
$\xi' = \xi \circ h$ et $g' = g \circ h$.
Le foncteur associant \`a
$Y$ de $\RDef_{Z}$ l'objet  $(Y, 0, 0)$ de $\RDefe_{Z}$, avec $0$ le morphisme
de valeur constante  $0$ dans $h[0,1,0]$, resp. $h[1, 0,0]$,
\'etant pleinement fid\`ele, on peut consid\'erer  $\RDef_{Z}$ comme une sous-cat\'egorie pleine de 
$\RDefe_{Z}$.
On associe \`a la
cat\'egorie $\RDefe_{Z}$ un anneau de Grothendieck
$K_0 (\RDefe_{Z})$ d\'efini de la fa\c con suivante.
Comme groupe ab\'elien c'est le
quotient du groupe ab\'elien libre sur les symboles $[Y \rightarrow Z,
\xi,g]$ pour $(Y \rightarrow Z, \xi,g)$ dans $\RDefe_{Z}$ par les
quatre
 relations suivantes :
\begin{equation}\label{r1}
[Y \rightarrow Z, \xi,g] = [Y' \rightarrow Z, \xi',g']
\end{equation}
si  $(Y \rightarrow Z, \xi,g)$ est isomorphe \`a $(Y' \rightarrow Z,
\xi',g')$,
\begin{equation}\label{r2}
[(Y \cup Y') \rightarrow Z, \xi, g]
+ [(Y \cap Y') \rightarrow Z,
\xi_{|Y \cap Y'},g_{|Y \cap Y'}] 
= [Y \rightarrow Z, \xi_{|Y},
g_{|Y}] + [Y' \rightarrow Z, \xi_{|Y'}, g_{|Y'}]
\end{equation}
si $Y$ et  $Y'$ sont des sous-assignements d\'efinissables d'un m\^eme $X$ de
$\RDef_Z$
et  $\xi$, $g$ sont d\'efinis sur  $Y \cup Y'$,
 \begin{equation}\label{r3}
 [Y\to Z,\xi,g+h]=[Y\to Z,\xi + \overline h,g]
 \end{equation}
pour $h:Y\to h[1,0,0]$ un morphisme d\'efinissable tel que  $\ord(h(y))  \geq 0$ pour tout
$y$ dans $Y$ et $\overline h$ la r\'eduction de $h$ modulo $t$,
et
 \begin{equation}\label{r4}
 [Y[0,1,0]\to Z,\xi+p,g]=0
 \end{equation}
avec $p:Y[0,1,0]\to h[0,1,0]$ la projection  et 
$Y[0,1,0]\to Z$, $g$, et  $\xi$ qui se factorisent \`a travers la projection
$Y[0,1,0]\to Y$.

On munit  $K_0 (\RDefe_{Z})$ d'une structure d'anneau en posant
$
[Y\to Z, \xi,g] \cdot  [Y'\to Z, \xi',g'] = [Y \otimes_Z Y'\to Z,
\xi\circ p_Y + \xi'\circ p_{Y'}, g\circ p_Y + g'\circ p_{Y'}]$,
avec $Y \otimes_Z Y'$ le produit fibr\' e de $Y$ et de $Y'$, $p_Y$
la projection sur  $Y$, et $p_{Y'}$ la projection sur $Y'$.

On \'ecrit  $\ee^\xi E(g) [Y\to Z]$ pour
$[Y\to Z, \xi,g]$. On abr\`ege $\ee^{0}E(g)[Y\to Z]$,
resp.~$\ee^\xi E(0)[Y\to Z]$, et $\ee^{0}E(0)[Y\to Z]$, en
$E(g)[Y\to Z]$, resp.~$\ee^\xi[Y\to Z]$, et  $[Y\to Z]$. De m\^eme on \'ecrit $\ee^\xi E(g)$ pour $\ee^\xi E(g)[Z\to Z]$, $E(g)$ pour $\ee^0 E(g)[Z\to Z]$ et $\ee^\xi$ pour $\ee^\xi E(0)[Z\to Z]$. L'\'el\'ement  $[Z\to Z]$
est l'unit\'e multiplicative de $K_0
(\RDefe_{Z})$.
On a une injection canonique
$K_0 (\RDef_Z ) \rightarrow
K_0 (\RDefe_Z )$ envoyant $[Y\to Z]$ sur  $[Y\to Z]$.
Pour $f : Z \rightarrow Z'$ dans  $\Def_{k}$, le produit fibr\' e induit un morphisme canonique
 $f^* : K_0 (\RDefe_{Z'}) \rightarrow
K_0 (\RDefe_{Z})$.
Si $f$ est un morphisme dans   $\RDef_{Z'}$,
la composition avec $f$ induit  un morphisme $f_! : K_0
(\RDefe_{Z}) \rightarrow K_0 (\RDefe_{Z'})$.

\subsection{Fonctions constructibles exponentielles}\label{cstre}Pour
$Z$ dans $\Def_k$ on d\'efinit l'anneau $\cC  (Z)^{\rm exp}$ des
fonctions constructibles exponentielles par
$\cC  (Z)^{\rm exp}:=\cC (Z)\otimes_{K_0 (\RDef_{Z})} K_0 (\RDefe_{Z})$. 
Notons que l'\'el\'ement
$E({\rm id})$ de  $\cC  (h[1,0,0])^{\rm exp}$ peut \^etre vu comme un caract\`ere additif
sur le corps valu\'e, non trivial par la relation (\ref{r3}).
De m\^eme,
l'\'el\'ement  $\ee^{\rm id}$ de $\cC (h[0,1,0])^{\rm exp}$
peut \^etre vu comme un caract\`ere additif du corps r\'esiduel, non trivial par la relation (\ref{r4}).
Pour tout $d\geq 0$ on d\'efinit  $\cC ^{\leq d} (Z)^{\rm exp} $ comme l'id\'eal de $\cC  (Z)^{\rm exp}$ engendr\'e par les fonctions caract\'eristiques $\11_{Z'}$ des sous-assignements d\'efinissables $Z'\subset Z$ de dimension
$\leq d$.
On pose $C (Z)^{\rm exp} := \oplus_d  C^d  (Z)^{\rm exp}$ avec
$C^d  (Z)^{\rm exp} := \cC ^{\leq d} (Z)^{\rm exp} / \cC ^{\leq
d-1} (Z)^{\rm exp}$.
C'est un groupe ab\'elien gradu\'e et un  $\cC  (Z)^{\rm
exp}$-module. Les \'el\'ements de $C (Z)^{\rm exp}$ sont les Fonctions
constructibles exponentielles.
Pour $S$ dans $\Def_k$ et $Z$ dans  $\Def_S$ on d\'efinit le groupe $ {\rm I}_S C (Z)^{\rm exp}$ des Fonctions
constructibles exponentielles $S$-int\'egrables comme 
 ${\rm I}_S C  (Z)^{\rm exp}:={\rm I}_S C
 (Z)\otimes_{K_0 (\RDef_{Z})} K_0 (\RDefe_{Z})$. C'est un sous-groupe gradu\'e de
$C  (Z)^{\rm exp}$.
On v\'erifie que les morphismes naturels $\cC (Z)\to \cC (Z)^{\rm exp}$, $\cC ^{\leq
d} (Z) \rightarrow \cC ^{\leq d} (Z)^{\rm exp}$, $C (Z)
\rightarrow C (Z)^{\rm exp}$ et ${\rm I}_S C (Z)\rightarrow  {\rm I}_S C
(Z)^{\rm exp}$ sont injectifs.
\subsection{Op\'erations naturelles}\label{opnat}
Soit $f : Z \rightarrow Z'$ un morphisme dans $\Def_k$. Les morphismes $f^*$
sur  $K_0 (\RDefe_{Z'})$ et sur  $\cC(Z')$ \'etant compatibles, on en d\'eduit par produit tensoriel
un morphisme naturel $f^* : \cC (Z')^{\rm
exp} \rightarrow \cC (Z)^{\rm exp}$.
Soit $Y$ dans $\Def_k$ et soit $Z$ un sous-assignement d\'efinissable de $Y[0,r,0]$, pour un 
$r\geq 0$. Soit
$f : Z \rightarrow Y$  le morphisme induit par la projection.
Le morphisme  $f_!: K_0 (\RDefe_{Z})
\rightarrow K_0 (\RDefe_{Y})$ induit un morphisme
d'anneaux
$f_! : \cC (Z)^{\rm exp} \rightarrow \cC (Y)^{\rm exp}$, 
ainsi qu'un morphisme de groupes
$f_! : C (Z)^{\rm exp} \rightarrow C
(Y)^{\rm exp}$. Si $Y$ est dans  $\Def_S$, le morphisme $f_!$ induit un morphisme
 $f_! :{\rm I}_S C (Z)^{\rm exp} \rightarrow {\rm I}_S C
(Y)^{\rm exp}$.

\section{L'\'enonc\'e principal}\label{sec7e}

\begin{theoreme}\label{mte}Soit $S$ dans
$\Def_k$. Il existe un unique foncteur de la cat\'egorie  $\Def_S$ vers celle des groupes ab\'eliens,
associant \`a un objet $Z$ le groupe ${\rm I}_S C
(Z)^{\rm exp}$, \`a un morphisme $f : Z \rightarrow Y$ in
$\Def_S$ le morphisme $f_!  : {\rm I}_S C  (Z)^{\rm exp} \rightarrow
{\rm I}_S C  (Y)^{\rm exp}$ et satisfaisant les axiomes suivants 

\begin{enumerate}

\item[]{{\rm A1 (Compatibilit\'e): }}Pour tout morphisme $f : Z
\rightarrow Y$ dans  $\Def_S$, le morphisme $f_!: {\rm I}_S C (Z)^{\rm
exp} \rightarrow {\rm I}_S C  (Y)^{\rm exp}$ est compatible avec les inclusions de groupes
 ${\rm I}_S C  (Z) \to {\rm I}_S C
(Z)^{\rm exp}$ et  ${\rm I}_S C  (Y) \to {\rm I}_S C (Y)^{\rm
exp}$ et le morphisme $f_!:{\rm I}_S C  (Z)\to {\rm I}_S C  (Y)$
construit dans  \cite{cl}.

\item[]{{\rm A2 (R\'eunion disjointe): }}Consid\'erons  $Z$ et  $Y$ dans $\Def_S$ et supposons que
 $Z$, resp. $Y$, est la r\'eunion disjointe de deux sous-assignements d\'efinissables $Z_1$ et $Z_2$, resp. $Y_1$ et $Y_2$,
de  $Z$, resp. $Y$. Pour tout morphisme $f : Z \rightarrow Y$ de $\Def_S$, avec $f (Z_i) \subset Y_i$ pour
$i = 1, 2$, de restriction $f_i : Z_i \to Y_i$, 
on a  $f_! = f_{1
!} \oplus f_{2!}$, modulo l'isomorphisme
${\rm I}_S C  (Z)^{\rm exp} \simeq {\rm I}_S C (Z_1)^{\rm exp}
\oplus {\rm I}_S C  (Z_2)^{\rm exp}$ et celui correspondant pour $Y$.

\item[]{{\rm A3 (Formule de projection): }}Pour tout morphisme $f : Z
\rightarrow Y$ dans $\Def_S$ et tout $\alpha$ dans $\cC  (Y)^{\rm
exp}$ et $\beta$ dans ${\rm I}_S C  (Z)^{\rm exp}$, 
si $f^*(\alpha) \beta$ est dans  ${\rm I}_S C  (Z)^{\rm exp}$,
alors  $f_! (f^*(\alpha) \beta) = \alpha
f_! (\beta)$.

\item[]{{\rm A4 (Projection le long des variables r\'esiduelles): }}  Supposons que
$f$ est la projection $f : Z = Y [0, n, 0] \rightarrow Y$ pour $Y$ dans
$\Def_S$ et que $\varphi$ est dans  ${\rm I}_S C  (Z)^{\rm
exp}$. Alors, $f_!(\varphi)$ est donn\'e par la construction de \ref{opnat}.

\item[]{{\rm A5 (Boules relatives de grand volume): }}Consid\'erons  $Y$ dans
$\Def_S$ et des morphismes d\'efinissables $\alpha : Y \rightarrow
\ZZ$ et  $\xi : Y \rightarrow h_{\GG_{m, k}}$, avec $\GG_{m, k}$ le groupe multiplicatif
 $\AA^1_k \setminus \{0\}$. 
Supposons que  $[\11_Z]$ est dans 
${\rm I}_S C (Z)^{\rm exp}$ et  que $Z$ est le sous-assignement d\'efinissable de $Y [1, 0, 0]$ d\'efini par
$\ord (z) =
\alpha (y)$ et $\ac (z) = \xi (y)$, et que  $f : Z \rightarrow
Y$ est le morphisme induit par la projection $Y[1,0,0] \rightarrow
Y$. Si, de plus, $\alpha(y)<0$ pour tout  $y$ dans $Y$, alors $f_! (E(z)[\11_Z]) = 0$.
\end{enumerate}
\end{theoreme}

\subsection*{Principe de la preuve:}Comme pour la preuve du th\'eor\`eme principal de \cite{cl},
le point-cl\'e
est la construction de $f_!$ lorsque $f$ est une projection
$Y [r, 0, 0]Ê\rightarrow Y$. Quand  $r = 1$, on utilise le th\'eor\`eme de d\'ecomposition cellulaire de Denef-Pas, cf. \cite{dp} et \cite{cl}, pour construire $f_!$ et on v\'erifie que la construction  est ind\'ependante de la d\'ecomposition cellulaire consid\'er\'ee.
Pour d\'eduire le cas $r > 1$ du cas $r = 1$, on \'etablit  une forme
convenable du th\'eor\`eme de Fubini  ainsi qu'un th\'eor\`eme de changement de variables pour $r = 1$.

Le th\'eor\`eme \ref{mte}  se g\'en\'eralise au cas relatif de la fa\c con suivante.
On fixe $\Lambda$ dans $\Def_k$ qui joue le r\^ole d'un espace de param\`etres.
Pour $Z$ dans $\Def_{\Lambda}$,
on consid\`ere, de fa\c con analogue \`a \cite{cr2},
l'id\'eal $\cC^{\leq d} (Z \rightarrow \Lambda)^{\rm exp}$
de $\cC (Z)^{\rm exp}$ engendr\'e
par les fonctions $\11_{Z'}$ avec $Z'$ sous-assignement d\'efinissable de $Z$
tel que toutes les fibres de $Z' \rightarrow \Lambda$ soient de dimension $\leq d$.
On pose $C (Z \rightarrow \Lambda)^{\rm exp} = \oplus_d  C^d (Z \rightarrow \Lambda)^{\rm exp}$ avec
$C^d (Z \rightarrow \Lambda) := \cC^{\leq d} (Z \rightarrow \Lambda)  / \cC^{\leq d-1} (Z \rightarrow \Lambda)$. Ce groupe ab\'elien gradu\'e s'identifie naturellement \`a
$C (Z \rightarrow \Lambda) \otimes_{K_0 (\RDef_Z)} K_0 (\RDef_Z^{\rm exp})$.
Pour $Z \rightarrow S$ un morphisme dans
$\Def_{\Lambda}$, on
 pose ${\rm I}_S C (Z \rightarrow \Lambda)^{\rm exp} := {\rm I}_S C (Z \rightarrow \Lambda)
\otimes_{K_0 (\RDef_Z)} K_0 (\RDef_Z^{\rm exp})$.

\begin{theoreme}\label{mter}Soit $\Lambda$ dans $\Def_k$ et soit $S$ dans $\Def_{\Lambda}$.
Il existe un unique foncteur de $\Def_S$ dans la cat\'egorie des
groupes ab\'eliens,
associant \`a tout morphisme $f : Z \rightarrow Y$ dans
$\Def_S$ un morphisme $f_{! \Lambda} :
{\rm I}_S C (Z\rightarrow \Lambda)^{\rm exp}
\rightarrow {\rm I}_S C (Y\rightarrow \Lambda)^{\rm exp}$
v\'erifiant les analogues de \textup{(A1)-(A5)} obtenus
en rempla\c cant
$ C (\_)$ par $C (\_ \rightarrow \Lambda)$.
\end{theoreme}

Soit $p : Z \rightarrow \Lambda$ un morphisme dans $\Def_k$, dont on suppose que toutes les fibres sont de
dimension $d$.
On note
$\cI_{\Lambda} (Z)^{\rm exp}$, ou encore $\cI_{p} (Z)^{\rm exp}$,
le  $\cC (\Lambda)^{\rm exp}$-module des fonctions  $\varphi$ in $\cC (Z)^{\rm exp}$ dont la classe $[\varphi]$
dans
$C^d (Z \rightarrow \Lambda)^{\rm exp}$ est int\'egrable  rel. $\Lambda$, i.e. est
dans ${\rm I}_{\Lambda} C ( Z \rightarrow \Lambda)^{\rm exp}$.
On \'ecrira 
$\mu_{\Lambda} (\varphi)$ ou 
$\mu_{p} (\varphi)$ pour
$p_{! \Lambda} ([\varphi])$.

\section{Transformation de Fourier}\label{fourier}

\subsection{Transformation de Fourier sur les variables valu\'ees}On fixe  $\Lambda$ dans $\Def_k$ et  un entier $d \geq 0$.
On consid\`ere le sous-assignement  d\'efinissable $\Lambda [2d, 0, 0]$, les
$d$ premi\`eres coordonn\'es sur les variables valu\'ees \'etant not\'ees
 $x = (x_1, \dots, x_d)$ et les $d$ derni\`eres
 $y = (y_1, \dots, y_d)$,
et on note $p_1 : \Lambda [2d, 0, 0] \rightarrow \Lambda [d, 0, 0]$ et
$p_2 : \Lambda [2d, 0, 0] \rightarrow \Lambda [d, 0, 0]$
la projection   sur les variables   $x$ et  $y$ respectivement.
On consid\`ere la fonction 
$E (xy) := E (\sum_{1 \leq i \leq d}Êx_i y_i)$ dans
$\cC (\Lambda [2d, 0, 0])^{\rm exp}$.
Si  $\varphi$ est une fonction dans 
$\cI_{\Lambda} (\Lambda [d, 0, 0])^{\rm exp}$, 
la classe $[E (xy) p_2^* (\varphi)]$
de
$E (xy) p_2^* (\varphi)$ dans $C^d (p_1 : \Lambda [2d, 0, 0] \rightarrow \Lambda [d, 0, 0])^{\rm exp}$
est int\'egrable rel. $p_1$, ce qui  permet de
d\'efinir la transformation de Fourier
$
\Fv : \cI_{\Lambda} (\Lambda [d, 0, 0])^{\rm exp}
\rightarrow 
\cC (\Lambda [d, 0, 0])^{\rm exp}
$
par
$
\Fv (\varphi) := \mu_{p_1} (E (xy) p_2^* (\varphi))$,
pour
$\varphi$ 
dans $\cI_{\Lambda} (\Lambda [d, 0, 0])^{\rm exp}$.
L'op\'erateur 
$\Fv$ est $\cC (\Lambda)^{\rm exp}$-lin\'eaire.

\subsection{Convolution}
On note  $x + y$ le morphisme
$\Lambda [2d, 0, 0] \rightarrow \Lambda [d, 0, 0]$
donn\'e par  $(x_1, \dots, x_d, y_1, \dots, y_d) \mapsto
(x_1 + y_1, \dots, x_d + y_d)$.

\begin{proposition}Soient  $f$ et  $g$ deux fonctions dans
$\cI_{\Lambda} (\Lambda [d, 0, 0])^{\rm exp}$.
\begin{enumerate}
\item[(1)]La fonction
$p_1^* (f) p_2^* (g)$ 
appartient \`a
$\cI_{x + y} (\Lambda [2d, 0, 0])^{\rm exp}$ et la fonction
$$
f \ast g := \mu_{x + y}Ê(p_1^* (f) p_2^* (g))
$$
appartient \`a
$\cI_{\Lambda} (\Lambda [d, 0, 0])^{\rm exp}$.
\item[(2)]Le produit de convolution
$(f, g) \mapsto f \ast g$ est  $\cC (\Lambda)^{\rm exp}$-lin\'eaire et munit 
$\cI_{\Lambda} (\Lambda [d, 0, 0])^{\rm exp}$
d'un produit associatif et commutatif.
\item[(3)]On a
$\Fv (f \ast g) = \Fv (f) \, \Fv (g)$.
\end{enumerate}
\end{proposition}

Pour $\varphi$ dans
$\cC (\Lambda [d, 0, 0])^{\rm exp}$ on \'ecrit
$\check{\varphi}$ pour $i^* (\varphi)$, avec  $i : \Lambda [d, 0, 0]
\rightarrow \Lambda [d, 0, 0]$ le $\Lambda$-morphisme envoyant
$x$ sur $- x$. Pour $\alpha : \Lambda \to \ZZ$ d\'efinissable, on note
$\varphi_{\alpha}$ la fonction caract\'eristique de l'ensemble des
$(\lambda, (x_i))$ avec $\ord \, x_i \geq \alpha (\lambda)$ pour tout $i$.

On a l'\'enonc\'e g\'en\'eral suivant (inversion de Fourier partielle) :

\begin{proposition}\label{parfouinv}
Soit  $\varphi$ une fonction dans  $\cI_{\Lambda} (\Lambda [d, 0, 0])^{\rm exp}$.
Pour tout  $\alpha$ d\'efinissable,
$\varphi_{\alpha} \Fv (\varphi)$ appartient \`a  $\cI_{\Lambda} (\Lambda [d, 0, 0])^{\rm exp}$
et 
$
\Fv (\varphi_{\alpha} \Fv (\varphi))=
\LL^{- \alpha d} \, \check{\varphi} \ast  \varphi_{- \alpha + 1}
$.
\end{proposition}

\subsection{Fonctions de Schwartz-Bruhat}

On d\'efinit l'espace
$\cS_{\Lambda} (\Lambda [d, 0, 0])^{\rm exp}$ des fonctions de
Schwartz-Bruhat sur  $\Lambda$
comme le  $\cC (\Lambda)^{\rm exp}$-sous-module
de $\cI_{\Lambda}  (\Lambda [d, 0, 0])^{\rm exp}$ form\'e des fonctions  $f$
telles qu'il existe 
$\alpha_0 : \Lambda \to \ZZ$ d\'efinissable
telle que 
$f \cdot \varphi_{\alpha} = f$ pour $\alpha \leq - \alpha_0$
et
$
f \ast \varphi_{\alpha} = \LL^{- \alpha d}Ê\, f$ pour $\alpha \geq \alpha_0$.

On a l'\'enonc\'e d'inversion suivant pour les fonctions de Schwartz-Bruhat : 

\begin{theoreme}La transformation de Fourier induit un isomorphisme
$$
\Fv : \cS_{\Lambda} (\Lambda [d, 0, 0])^{\rm exp} \simeq \cS_{\Lambda} (\Lambda [d, 0, 0])^{\rm exp}
$$
tel que, pour toute fonction $\varphi$ dans
$\cS_{\Lambda} (\Lambda [d, 0, 0])^{\rm exp}$,
$
(\Fv \circ \Fv) (\varphi)
= \LL^{-d}\,  \check{\varphi}$.
\end{theoreme}

On en d\'eduit l'\'enonc\'e d'inversion g\'en\'eral suivant pour les fonctions int\'egrables \`a transform\'ee de Fourier int\'egrable : 

\begin{theoreme}Soit $\varphi$ dans
$\cI_{\Lambda} (\Lambda [d, 0, 0])^{\rm exp}$. On suppose que $\Fv (\varphi)$
est aussi dans 
$\cI_{\Lambda} (\Lambda [d, 0, 0])^{\rm exp}$. Alors les fonctions
$(\Fv \circ \Fv) (\varphi)$ et $\LL^{-d}\,  \check{\varphi}$ ont m\^eme classe
dans $C^d (\Lambda [d, 0, 0] \rightarrow \Lambda)$.
\end{theoreme}

\section{Sp\'ecialisation et transfert}

\subsection{Sp\'ecialisation}\label{nd}

On suppose d\'esormais que  $k$ est un corps de nombres d'anneaux des entiers $\cO$.
On note 
$\cA_{\cO}$ l'ensemble des compl\'etions $p$-adiques de
toutes les extensions finies de $k$ et
$\cB_{\cO}$ l'ensemble des corps locaux de caract\'eristique $>0$ qui sont des $\cO$-modules. 
Pour $K$ dans  $\cC_{\cO}:=\cA_{\cO}\cup \cB_{\cO}$, on d\'esigne par $R_K$ son anneau de valuation, $M_K$
l'id\'eal maximal, $k_K$  le corps r\'esiduel,
$q_K$ le cardinal de $k_K$, et $\varpi_K$ une uniformisante de $R_K$. Il existe un unique morphisme $\ac:K^\times \to
k_K^\times$ prolongeant $R_K^\times\to
k_K^\times$ et envoyant $\varpi_K$ sur $1$, que l'on \'etend par
$\ac(0)=0$. On note $\cD_K$ l'ensemble des caract\`eres additifs
$\psi:K\to\CC^\times$ tels que
  $\psi(x)=\exp((2\pi i /p)
{\rm Tr}_{k_K}
(\bar x))$ pour $x\in R_K$, avec $p$ la caract\'eristique de $k_K$,  
${\rm Tr}_{k_K}$
la trace de $k_K$ relativement \`a son sous-corps premier et  $\bar x$ la classe de $x$ dans $k_K$. Pour $N>0$, on d\'efinit $\cA_{\cO, N}$ comme l'ensemble des corps $K$ dans $\cA_{\cO}$ tels que
$k_K$ soit de caract\'eristique $> N$, et de m\^eme pour $\cB_{\cO, N}$ et $\cC_{\cO, N}$.

Pour pouvoir interpr\'eter nos formules dans les corps de
$\cC_\cO$, on restreint le langage $\LPas$ au sous-langage 
$\LO$
avec coefficients dans $k$ pour la sorte r\'esiduelle et dans
$\cO \llb t \rrb$ pour la sorte valu\'ee. On note $\Def(\LO)$ la sous-cat\'egorie de $\Def_k$
des objets d\'efinissables
 dans le langage
$\LO$. On proc\`ede de m\^eme pour les fonctions, etc., ainsi,
pour $S$ dans $\Def(\LO)$, on note $\cC (S, \LO)^{\rm exp}$  l'anneau des fonctions constructibles exponentielles d\'efinissables
dans
$\LO$.
On voit $K$ de $\cC_{\cO}$ comme une  $\cO \llb t \rrb$-alg\`ebre via le morphisme
$\lambda_{\cO,K}: \sum_{i\in\NN}a_it^i\mapsto
\sum_{i\in\NN}a_i\varpi_K^i,$
ainsi, si on interpr\`ete $a$ dans  $\cO \llb t \rrb $ par
$\lambda_{\cO,K}(a)$, toute  $\LO$-formule $\varphi$
definit pour  $K\in\cC_\cO$ un sous-ensemble d\'efinissable $\varphi_K$ d'un
$K^m \times k_K^n \times \ZZ^r$.
Il r\'esulte de r\'esultats 
d'Ax, Kochen et Er{\v s}ov
que si deux $\LO$-formules  $\varphi$ et $\varphi'$ d\'efinissent le m\^eme sous-assignement $X$ de
$h[m,n,r]$, alors $\varphi_K=\varphi'_{K}$ pour tout
$K$ dans  $\cC_{\cO,N}$ avec $N$ assez grand.
Avec abus de notation on note $X_K$ l'ensemble d\'efini
par $\varphi_K$
pour
$K$ dans  $\cC_{\cO,N}$ avec $N \gg 0$. De m\^eme, tout morphisme $\LO$-definissable
$f:X\to Y$
se sp\'ecialise en   $f_K:X_K\to Y_K$ pour $K$ dans
$\cC_{\cO,N}$ avec $N \gg 0$.
De fa\c con similaire,  $\varphi$ dans
$\cC(X,\LO)$ se sp\'ecialise en $\varphi_K : X_{K}\to \QQ$ pour $K$ dans
$\cC_{\cO,N}$ avec $N \gg 0$.
En effet consid\'erons
$\varphi$ dans $K_0(\RDef_X (\LO))$ de la forme  $[\pi:W\to X]$ avec $W$ dans
$\RDef_X(\LO)$. Pour  $K$ dans $\cC_{\cO, N}$, avec $N \gg 0$,
 on dispose de  $\pi_K:W_K\to X_K$ et on d\'efinit
 $\varphi_K:X_K\to\QQ$ par $: x\mapsto {\rm card} \left(\pi_K^{-1}(x)\right)$.
Pour  $\varphi$ dans  $\cP(X)$, on sp\'ecialise $\LL$ en
$q_K$ et  $\alpha:X\to\ZZ$ en
$\alpha_K:X_K\to\ZZ$. Par produit tensoriel on en d\'eduit la sp\'ecialisation
$\varphi \mapsto \varphi_K$ pour 
$\varphi$ dans $\cC(X,\LO)$. Pour le cas exponentiel, soit 
$\varphi$ dans $K_0(\RDef_X (\LO))^{\rm exp}$
de la forme $[W,g,\xi]$. Pour $\psi_K$ dans $\cD_K$, on sp\'ecialise
$\varphi$ en 
$\varphi_{K,\psi_K}:X_K\to\CC$ donn\'ee par $x
\mapsto \sum_{y\in \pi_K^{-1}(x)}\psi_K(g_K(y))\exp((2\pi i /p)
{\rm Tr}_{k_K}
(\xi_K(y)))$ pour $K$ dans
$\cC_{\cO,N}$ avec $N \gg 0$.
Par produit tensoriel on en d\'eduit la sp\'ecialisation
$\varphi \mapsto \varphi_K$ pour 
$\varphi$ dans $\cC(X,\LO)^{\rm exp}$.

Soit  $K$ dans  $\cC_\cO$ et soit $A$ une partie de $K^m\times
k_K^n \times \ZZ^r$. On d\'efinit la dimension de $A$ comme la dimension de l'adh\'erence de Zariski $\bar A$ de la projection de $A$ dans $\AA_K^m$. Pour $d \geqÔ{\rm dim} A$,
on obtient une mesure $\mu_d$ sur $A$ (nulle si $d >Ô{\rm dim} A$)
par restriction du produit de la mesure $d$-dimensionnelle canonique sur $\bar A (K)$ avec la mesure de comptage sur $k_K^n \times \ZZ^r$.

Fixons un morphisme 
$f : S \rightarrow \Lambda$ dans $\Def (\LO)$. 
Consid\'erons $\phi$ dans
$\cC^{\leq d} (S \rightarrow \Lambda, \LO)^{\rm exp}$.
On d\'emontre que pour
$N \gg 0$, pour $K$ dans $\cC_{\cO,N}$, et $\psi_K$ dans $ \cD_K$,
pour tout $\lambda$ dans $\Lambda_K$,
le support 
de la restriction
$\phi_{K,\psi_K, \lambda}$ 
de
$\phi_{K,\psi_K}$
\`a $f_K^{-1}Ê(\lambda)$ 
est de dimension $\leq d$ et que de plus si la classe $\varphi$ de $\phi$ dans 
$C^{d} (S \rightarrow \Lambda, \LO)^{\rm exp}$ est relativement int\'egrable, les 
$\phi_{K,\psi_K,
\lambda}$ sont toutes $\mu_d$-int\'egrables.
On note alors
$\mu_{\Lambda_K}(\varphi_{K,\psi_K})$ la fonction sur $\Lambda_K$ d\'efinie par 
$\lambda \mapsto \mu_d (\phi_{K,\psi_K, \lambda})$, qui ne d\'epend que de $\varphi$ pour $N \gg 0$.
Cette construction s'\'etend par lin\'earit\'e \`a
${\rm I} C (S \rightarrow \Lambda, \LO)^{\rm exp}$.

\begin{theoreme}[Principe de sp\'ecialisation]\label{compres3}Soit
$f : S \rightarrow \Lambda$ un morphisme dans $\Def (\LO)$. Soit 
$\varphi$ dans  
${\rm I}_{\Lambda} C (S \rightarrow \Lambda, \LO)^{\rm exp}$. Il existe
$N>0$ tel que, pour  tout $K$ dans $\cC_{\cO,N}$ et tout $\psi_K$ dans $\cD_K$, $\left(\mu_\Lambda(\varphi)\right)_{K,\psi_K}
=
\mu_{\Lambda_K}(\varphi_{K,\psi_K})$.
\end{theoreme}

\subsection{Principe de transfert pour les int\'egrales avec param\`etres}

\begin{theoreme}\label{strong} Soit $\varphi$ dans
$\cC (\Lambda, \LO)^{\rm exp}$. Il existe un entier  $N$
tel que pour tous $K_1$, $K_2$ dans  $\cC_{\cO,N}$ avec $k_{K_1}\simeq
k_{K_2}$, 
$$\varphi_{K_1,\psi_{K_1}} = 0 \quad \mbox{pour tout} \quad 
\psi_{K_1}\in\cD_{K_1}
$$
si et seulement si
$$ \varphi_{K_2,\psi_{K_2}} = 0  \quad \mbox{pour tout} \quad 
\psi_{K_2}\in\cD_{K_2}.
$$
\end{theoreme}

On d\'eduit des th\'eor\`emes \ref {compres3} et \ref{strong}:

\begin{theoreme}[Principe de transfert pour les int\'egrales avec param\`etres]\label{strongaxkc}
Soient $S \rightarrow \Lambda$ et  $S' \rightarrow \Lambda$ des morphismes dans
 $\Def (\LO)$. Soit  $\varphi$ dans  ${\rm I}_{\Lambda} C (S
\rightarrow \Lambda, \LO)^{\rm exp}$ et $\varphi'$ dans ${\rm I}_{\Lambda} C
(S' \rightarrow \Lambda, \LO)^{\rm exp}$. Il existe un entier $N$ tel que pour tous les corps
 $K_1$ et $K_2$ dans  $\cC_{\cO,N}$ avec
$k_{K_1}\simeq k_{K_2}$, 
$$
\mu_{ \Lambda_{K_1}} (\varphi_{K_1,\psi_{K_1}}) = \mu_{
\Lambda_{K_1}} (\varphi'_{K_1,\psi_{K_1}})  \quad \mbox{pour tout} \quad 
\psi_{K_1}\in\cD_{K_1}
 $$
si et seulement si
 $$
\mu_{ \Lambda_{K_2}} (\varphi_{K_2,\psi_{K_2}}) = \mu_{
\Lambda_{K_2}} (\varphi'_{K_2,\psi_{K_2}}) \quad \mbox{pour tout} \quad 
\psi_{K_2}\in\cD_{K_2}.
$$
\end{theoreme}

En l'absence d'exponentielles une forme du th\'eor\`eme
pr\'ec\'edent se trouve d\'ej\`a dans \cite{miami}.
Il devrait s'appliquer aux int\'egrales 
orbitales apparaissant dans diff\'erentes formes du Lemme Fondamental.
Notons que notre approche permet de s'affranchir de conditions de constance locale, telles celles figurant dans
 \cite{CH}. Rappelons que le Lemme Fondamental pour les groupes unitaires a \'et\'e d\'emontr\'e par
 Laumon et Ng\^o
 \cite{ln} sur les corps de fonctions et que Waldspurger en a d\'eduit le cas des corps $p$-adiques \cite{wal}.
 D'autre part notre r\'esultat s'applique aux int\'egrales
apparaisssant dans la conjecture de Jacquet-Ye \cite{JY},
d\'emontr\'ee par Ng\^o \cite{Ngo} sur les corps de fonctions  et par Jacquet \cite{J} en g\'en\'eral.

\subsection*{}\begin{small}Pendant la r\'ealisation de ce projet, le premier auteur \'etait chercheur postdoctoral du Fonds de Recherche Scientifique - Flandres (Belgique) et il a b\'en\'efici\'e de la  bourse Marie Curie de la Commission Europ\' eenne HPMF - CT 2005-007121. 
\end{small}

\bibliographystyle{amsplain}

\end{document}